\documentclass[11pt]{article}
\usepackage{latexsym,amssymb,amsmath,geometry,cite,enumerate}
\newtheorem{theorem}{Theorem}

\newtheorem{claim}{Claim}

\allowdisplaybreaks

\usepackage{lineno}
\usepackage{setspace}

\begin{document}

\onehalfspace

\title{Large Values of the Clustering Coefficient}

\author{Michael Gentner$^1$, Irene Heinrich$^2$, Simon J\"{a}ger$^1$, Dieter Rautenbach$^1$}

\date{}

\maketitle

\begin{center}
$^1$ Institute of Optimization and Operations Research, Ulm University, Ulm, Germany, 
\texttt{michael.gentner, simon.jaeger, dieter.rautenbach@uni-ulm.de}\\[3mm]
$^2$ Department of Mathematics, University of Kaiserslautern, Kaiserslautern, Germany,
\texttt{heinrich@mathematik.uni-kl.de}
\end{center}

\begin{abstract}
A prominent parameter in the context of network analysis, 
originally proposed by Watts and Strogatz (Collective dynamics of `small-world' networks, Nature 393 (1998) 440-442),
is the clustering coefficient of a graph $G$.
It is defined as the arithmetic mean of the clustering coefficients of its vertices, 
where the clustering coefficient of a vertex $u$ of $G$ is the relative density $m(G[N_G(u)])/{d_G(u)\choose 2}$ 
of its neighborhood if $d_G(u)$ is at least $2$, and $0$ otherwise. 
It is unknown which graphs maximize the clustering coefficient among all connected graphs of given order and size.

We determine the maximum clustering coefficients 
among all connected regular graphs of a given order,
as well as  
among all connected subcubic graphs of a given order.
In both cases, we characterize all extremal graphs.
Furthermore, we determine the maximum increase of the clustering coefficient caused by adding a single edge.
\end{abstract}

{\small 

\begin{tabular}{lp{13cm}}
{\bf Keywords:} & clustering coefficient; connected caveman graph; cliquishness
\end{tabular}
}

\pagebreak

\section{Introduction}

Watts and Strogatz \cite{wast} proposed the clustering coefficient of a graph in order to quantify the corresponding property of networks.
For a vertex $u$ of a graph $G$, let 
the {\it clustering coefficient} of $u$ in $G$ be
$$C_u(G)=
\begin{cases}
\frac{m(G[N_G(u)])}{{d_G(u)\choose 2}} &, \mbox{ if $d_G(u)\geq 2$,}\\
0 &, \mbox{ otherwise,}
\end{cases}
$$
where $m(G[N_G(u)])$ is the size of the subgraph of $G$ induced by the neighborhood $N_G(u)$ of $u$ in $G$,
that is, $m(G[N_G(u)])$ equals exactly the number of triangles of $G$ that contain $u$.
Furthermore, let the {\it clustering coefficient} of $G$ be the average 
$$C(G)=\frac{1}{n(G)}\sum\limits_{u\in V(G)}C_u(G)$$
of the clustering coefficients of its $n(G)$ vertices.

While the clustering coefficient received a lot of attention within social network analysis \cite{bojapr,scwa,wa1,wa2},
some fundamental mathematical problems related to it are still open.
It is unknown, for instance, which graphs maximize the clustering coefficient among all connected graphs of a given order and size.

Watts \cite{wa1,wa2} suggested the so-called {\it connected caveman} graphs as a possible extremal construction.
For integers $k$ and $\ell$ at least $2$,
these arise from $\ell$ disjoint copies $G_1,\ldots,G_\ell$ of $K_{k+1}-e$ arranged cyclically
by adding, for every $i$ in $[\ell]$, 
an edge between 
one of the two vertices of degree $k-1$ in $G_i$ and 
one of the $k-1$ vertices of degree $k$ in $G_{i+1}$, where the indices are identified modulo $\ell$.
Actually, it is rather obvious that these graphs do not have the largest clustering coefficient 
among all connected graphs of given order and size,
because removing the edge between $G_1$ and $G_2$,
and adding a new edge between the two vertices of degree $k-1$ in $G_1$,
increases the clustering coefficient.

Fukami and Takahashi \cite{futa1,futa2} considered {\it clustering coefficient locally maximizing graphs}
whose clustering coefficient cannot be increased by some local operations such as an edge swap.

In the present paper 
we determine the maximum clustering coefficients 
among all connected regular graphs of a given order,
as well as  
among all connected subcubic graphs of a given order.
In both cases, we characterize all extremal graphs.
Furthermore, we determine the maximum increase of the clustering coefficient caused by adding a single edge.

\section{Results}

We introduce a slightly modified version of the connected caveman graphs.
For integers $k$ and $\ell$ with $k\geq 3$ and $\ell\geq 2$,
let $G(k,\ell)$ be the $k$-regular connected graph that arises 
from $\ell$ disjoint copies $G_1,\ldots,G_\ell$ of $K_{k+1}-e$ arranged cyclically
by adding, for every $i$ in $[\ell]$, 
an edge between a vertex in $G_i$ and a vertex in $G_{i+1}$, 
where the indices are identified modulo $\ell$.
Note that $G(k,\ell)$ is uniquely determined up to isomorphism by the requirement of $k$-regularity.

\begin{theorem}\label{theorem1}
Let $k$ and $n$ be integers with $n\geq k+2$ and $k\geq 3$.
If $G$ is a connected $k$-regular graph of order $n$, then 
$$C(G)\leq 1-\frac{6}{k(k+1)}$$
with equality if and only if $n/(k+1)$ is an integer and $G$ equals $G(k,n/(k+1))$.
\end{theorem}
{\it Proof:} Let $G$ be a connected $k$-regular graph of order $n$.
For a non-negative integer $i$, let $V_i$ be the set of vertices $u$ of $G$ with $m(G[N_G(u)])={k\choose 2}-i$.
Since $G$ is connected and has order at least $k+2$, no vertex has a complete neighborhood, that is, $V_0$ is empty.
For a set $U$ of vertices of $G$, let $\sigma(U)=\sum\limits_{u\in U}C_u(G)$.
In order to obtain a useful decomposition of $G$, we consider some special graphs.
For $k\leq 4$, one such graph suffices, while for $k\geq 5$, two more are needed.

Let $G_1,\ldots,G_r$ be a maximal collection of disjoint subgraphs of $G$ that are all copies of $K_{k+1}-e$.
Let $A=V(G_1)\cup \cdots \cup V(G_r)$ and $R=V(G)\setminus A$.
Note that every vertex in $A$ has at most one neighbor in $R$.
Suppose that $R$ contains a vertex $u$ from $V_1$.
Since every vertex in $N_G(u)$ has at least two neighbors in the closed neighborhood $N_G[u]$ of $u$,
the subgraph $G_{r+1}$ of $G$ induced by $N_G[u]$ 
does not intersect $A$ 
and is a copy of $K_{k+1}-e$.
Now, $G_1,\ldots,G_r,G_{r+1}$ contradicts the maximality of the above collection,
which implies that $R$ does not intersect $V_1$.

Since each $G_i$ contains $k-1$ vertices from $V_1$ and two vertices whose neighborhood induces $K_1\cup K_{k-1}$,
and $|A|=r(k+1)$,
we have
\begin{eqnarray*}
\sigma(A)=\sum\limits_{i\in [r]}\sum\limits_{u\in V(G_i)}C_u(G) 
= r\left((k-1)\frac{{k\choose 2}-1}{{k\choose 2}}+2\frac{{k-1\choose 2}}{{k\choose 2}}\right)
= |A|\left(1-\frac{6}{k(k+1)}\right).
\end{eqnarray*}
Since $R$ does not intersect $V_1$, we have
\begin{eqnarray*}
\sigma(R)=\sum\limits_{u\in R}C_u(G) 
\leq |R|\frac{{k\choose 2}-2}{{k\choose 2}}
=|R|\left(1-\frac{4}{(k-1)k}\right).
\end{eqnarray*}
First, let $k\leq 4$.
Since $1-\frac{4}{(k-1)k}<1-\frac{6}{k(k+1)}$ in this case, we obtain
\begin{eqnarray*}
C(G) &=& \frac{1}{n(G)}(\sigma(A)+\sigma(R))\\
&\leq & \frac{1}{n(G)}(|A|+|R|)\left(1-\frac{6}{k(k+1)}\right)\\
&=& 1-\frac{6}{k(k+1)},
\end{eqnarray*}
with equality if and only if $A=V(G)$ and $R=\emptyset$, which implies that $k+1$ divides $n$, and $G$ equals $G(k,n/(k+1))$.

Now, let $k\geq 5$.
In this case, we need to refine the partition of $V(G)$ into $A$ and $R$ further.
Let $G_1,\ldots,G_r$, $A$, and $R$ be exactly as above, and recall that $R$ does not intersect $V_1$.
Let $H_1,\ldots,H_s$ be a maximal collection of disjoint subgraphs of $G[R]$
that are all copies of the two possible graphs 
that arise from $K_{k+1}$ by removing two edges.
Let $B=V(H_1)\cup \cdots \cup V(H_s)$, and $S=V(G)\setminus (A\cup B)$.
Note that every vertex in $A\cup B$ has at most two neighbors in $S$.
Suppose that $S$ contains a vertex $u$ from $V_2$.
Since $k\geq 5$, every vertex in $N_G(u)$ has at least three neighbors in the closed neighborhood $N_G[u]$ of $u$.
Hence, the subgraph $H_{s+1}$ of $G$ induced by $N_G[u]$ 
does not intersect $A\cup B$
and is a copy of one of the two possible graphs that arise from $K_{k+1}$ by removing two edges.
Now, $H_1,\ldots,H_s,H_{s+1}$ contradicts the maximality of the above collection,
which implies that $S$ does not intersect $V_1\cup V_2$.

If $H_i$ for some $i$ in $[s]$ is a copy of $K_{k+1}$ minus two non-incident edges, 
then $H_i$ contains $k-3$ vertices from $V_2$
and four vertices whose neighborhood induces a graph that arises from $K_1\cup (K_{k-1}-e)$ by adding at most two edges,
which implies 
\begin{eqnarray*}
\sigma(V(H_i)) 
& \leq & (k-3)\frac{{k\choose 2}-2}{{k\choose 2}}+4\frac{{k-1\choose 2}-1+2}{{k\choose 2}}
=\frac{k^3-13 k+28}{(k-1)k}.
\end{eqnarray*}
If $H_i$ for some $i$ in $[s]$ is a copy of $K_{k+1}$ minus two incident edges, 
then $H_i$ contains $k-2$ vertices from $V_2$,
two vertices whose neighborhood induces a graph that arises from $K_1\cup K_{k-1}$ by adding at most one edge, and 
one vertex whose neighborhood induces a graph that arises from $K_1\cup K_1\cup K_{k-2}$ by adding at most one edge,
which implies 
\begin{eqnarray*}
\sigma(V(H_i)) 
& \leq & (k-2)\frac{{k\choose 2}-2}{{k\choose 2}}+2\frac{{k-1\choose 2}+1}{{k\choose 2}}+\frac{{k-2\choose 2}+1}{{k\choose 2}}
=\frac{k^3-13 k+24}{(k-1)k}.
\end{eqnarray*}
Using $|B|=s(k+1)$, we obtain 
\begin{eqnarray*}
\sigma(B) & = & \sum\limits_{i\in [s]}\sigma(V(H_i))
\leq |B|\frac{k^3-13 k+28}{(k-1)k(k+1)}
= |B|\left(1-\frac{12k-28}{(k-1)k(k+1)}\right).
\end{eqnarray*}
Since $S$ does not intersect $V_1\cup V_2$, we have 
\begin{eqnarray*}
\sigma(S) & \leq & |S|\frac{{k\choose 2}-3}{{k\choose 2}}
= |S|\left(1-\frac{6}{(k-1)k}\right).
\end{eqnarray*}
Since 
$$\max\left\{ 1-\frac{12k-28}{(k-1)k(k+1)},1-\frac{6}{(k-1)k}\right\}<1-\frac{6}{k(k+1)}$$
for $k\geq 5$,
we obtain
\begin{eqnarray*}
C(G) &=& \frac{1}{n(G)}(\sigma(A)+\sigma(B)+\sigma(S))\\
&\leq &\frac{1}{n(G)}(|A|+|B|+|S|)\left(1-\frac{6}{k(k+1)}\right)\\
&=& 1-\frac{6}{k(k+1)},
\end{eqnarray*}
with equality if and only if $A=V(G)$ and $B=S=\emptyset$, 
which implies that $k+1$ divides $n$, and $G$ equals $G(k,n/(k+1))$.
This completes the proof. $\Box$

\medskip

\noindent Recall that the {\it diamond} is the unique graph with degree sequence $2,2,3,3$.

\begin{theorem}\label{theorem2}
If $G$ is a connected subcubic graph of order $n$ at least $6$, then 
\begin{eqnarray}\label{e1}
C(G)\leq 
\begin{cases}
\frac{7}{12}+\frac{12}{12n} & \mbox{, if $n\equiv 0$ mod $4$,}\\
\frac{7}{12}+\frac{13}{12n} & \mbox{, if $n\equiv 1$ mod $4$,}\\
\frac{7}{12}+\frac{14}{12n} & \mbox{, if $n\equiv 2$ mod $4$, and}\\
\frac{7}{12}+\frac{11}{12n} & \mbox{, if $n\equiv 3$ mod $4$.}\\
\end{cases}
\end{eqnarray}
\end{theorem}
{\it Proof:} Let ${\cal G}$ be the set of all connected subcubic graphs of order $n$ at least $6$. We assume that $G$ is chosen within ${\cal G}$ in such a way that 
\begin{enumerate}[(i)]
\item its clustering coefficient $C(G)$ is as large as possible,
\item subject to condition (i), the size $m(G)$ of $G$ is as small as possible, and,
\item subject to conditions (i) and (ii), the number of triangles in $G$ that contain at most one vertex of degree $2$ in $G$ is as small as possible.
\end{enumerate}
We establish a series of structural properties of $G$.

\begin{claim}\label{claim1}
Every subgraph $D$ of $G$ that is a diamond is induced and forms an endblock.
\end{claim}
{\it Proof of Claim \ref{claim1}:}
Since $G$ is subcubic, connected, and of order more than $4$, the subgraph $D$ is induced.
If all vertices in $D$ have degree $3$ in $G$,
then contracting $D$ to a single vertex $u$,
adding a new triangle $xyz$,
and adding the new edge $ux$
yields a graph $G'$ in ${\cal G}$ with $C(G')\geq C(G)+\frac{1}{3n}$,
contradicting the choice of $G$.
Note that the two neighbors of $u$ in $G'$ that are distinct from $x$ may be adjacent, 
in which case $C(G')>C(G)+\frac{1}{3n}$.
In view of the order, 
this implies that $D$ contains exactly one vertex of degree $2$,
and, hence, forms an endblock of $G$.
$\Box$

\begin{claim}\label{claim2}
Every edge of $G$ that lies in some cycle also lies in some triangle.
\end{claim}
{\it Proof of Claim \ref{claim2}:}
If the edge $uv$ of $G$ lies in some cycle but in no triangle,
then removing $uv$ yields a graph $G'$ in ${\cal G}$ with $C(G')\geq C(G)$ and $m(G')<m(G)$,
contradicting the choice of $G$.
Note that, if some triangle of $G$ contains $u$ or $v$, then $C(G')>C(G)$.
$\Box$

\begin{claim}\label{claim3}
Every block of $G$ is $K_2$, $K_3$, or a diamond.
\end{claim}
{\it Proof of Claim \ref{claim3}:}
Suppose, for a contradiction, 
that $B$ is a block of $G$ that is neither $K_2$, nor $K_3$, nor a diamond. 
Note that every edge of $B$ lies in some cycle,
and, hence, by Claim \ref{claim2},
also lies in some triangle.
Let $uvw$ be a triangle in $B$.
Since $B$ is not $K_3$,
we may assume that $u$ has a neighbor $x$ in $V(B)\setminus \{ v,w\}$.
Since the edge $ux$ of $B$ lies in some triangle, 
we may assume that $x$ and $v$ are adjacent.
Since $B$ is not a diamond, 
we may assume, by symmetry, that $x$ has a neighbor $y$ in $V(B)\setminus \{ u,v,w\}$.
Since the edge $xy$ of $B$ lies in some triangle,
we obtain that $y$ is adjacent to $u$ or $v$,
contradicting the assumption that $G$ is subcubic.
$\Box$

\medskip

\noindent Let ${\cal D}$ be the set of blocks of $G$ that are diamonds.
For $i\in \{ 2,3\}$, let ${\cal I}_i$ be the set of blocks of $G$ that are triangles that contain exactly $i$ vertices of degree $3$ in $G$.
Let ${\cal I}={\cal I}_2\cup {\cal I}_3$.
Finally, let ${\cal S}$ be the set of vertices of $G$ that do not lie in some triangle and have degree at most $2$ in $G$.

The triangles in ${\cal I}$ are called {\it inner triangles}.
By condition (iii), $G$ has as few inner triangles as possible 
given the other conditions.

If $u\in {\cal S}$ has degree $2$, then {\it resolving $u$} means to remove $u$ from $G$, and to connect its two neighbors by a new edge.
If $u\in {\cal S}$ has degree $1$, then {\it resolving $u$} simply means to remove $u$ from $G$.
Note that resolving some vertex from ${\cal S}$ yields a connected subcubic graph $G'$ of order $n-1$ with $C(G')\geq C(G)$.

\begin{claim}\label{claim4}
Either ${\cal D}$ or ${\cal I}$ is empty.
\end{claim}
{\it Proof of Claim \ref{claim4}:}
Suppose, for a contradiction,  
that $G$ contains a diamond $D$ and an inner triangle $T$.
By Claim \ref{claim1}, $D$ contains a vertex $u$ of degree $2$ in $G$.
Let $v$ be a neighbor of $u$.
If $T\in {\cal I}_2$, then 
contracting $T$ to a single vertex $w$, 
removing $u$, 
adding a new triangle $xyz$, and 
adding the new edge $xw$
yields a graph $G'$ in ${\cal G}$ with $C(G')=C(G)+\frac{1}{3n}$,
contradicting the choice of $G$.
If $T\in {\cal I}_3$, then 
contracting $T$ to a single vertex, 
removing $u$, 
adding a new triangle $xyz$, and 
adding the new edge $xv$
yields a graph $G'$ in ${\cal G}$ with $C(G')=C(G)+\frac{1}{3n}$,
contradicting the choice of $G$.
$\Box$

\begin{claim}\label{claim5}
$|{\cal D}|\leq 2$.
\end{claim}
{\it Proof of Claim \ref{claim5}:}
Suppose, for a contradiction, 
that $G$ has three blocks $B_1$, $B_2$, and $B_3$ that are diamonds.
By Claim \ref{claim1}, these are all endblocks, and
removing the three vertices, say $v_1$, $v_2$, and $v_3$, of degree $2$ in $G$ from $B_1$, $B_2$, and $B_3$,
adding a new triangle $xyz$, and
adding a new edge between $x$ and one of the two neighbors of $v_1$ in $B_1$
yields a graph $G'$ in ${\cal G}$ with $C(G')=C(G)+\frac{2}{3n}$,
contradicting the choice of $G$.
$\Box$

\begin{claim}\label{claim6}
${\cal S}$ is empty.
\end{claim}
{\it Proof of Claim \ref{claim6}:}
Suppose, for a contradiction,  
that ${\cal S}$ contains some vertex $u$.

If ${\cal I}_2$ contains a triangle $T$, then 
resolving $u$,
contracting $T$ to a single vertex $v$,
adding a new triangle $xyz$, and 
adding the new edge $xv$
yields a graph $G'$ in ${\cal G}$ with $C(G')\geq C(G)+\frac{2}{3n}$,
contradicting the choice of $G$.
If ${\cal I}_3$ contains a triangle $T$, then 
resolving $u$,
contracting $T$ to a single vertex $v$,
adding a new triangle $xyz$, and then
replacing one of the edges incident with $v$, say $vw$,
with the two new edges $xv$ and $yw$
yields a graph $G'$ in ${\cal G}$ with $C(G')\geq C(G)+\frac{2}{3n}$,
contradicting the choice of $G$.
Hence, we may assume that $G$ has no inner triangles.

If $G$ has a block $B$ that is a triangle, then $B$ is an endblock, and 
resolving $u$,
adding a new vertex $v$, and 
adding two new edges between $v$ and the two vertices of degree $2$ in $B$
yields a graph $G'$ in ${\cal G}$ with $C(G')\geq C(G)+\frac{1}{3n}$,
contradicting the choice of $G$.
Hence, we may assume that $G$ has no blocks that are triangles.

If $G$ has two blocks $B_1$ and $B_2$ that are diamonds, 
then $B_1$ and $B_2$ are endblocks by Claim \ref{claim1}, and
resolving $u$,
removing the two vertices, say $v_1$ and $v_2$, of degree $2$ in $G$ from $B_1$ and $B_2$, respectively,
adding a new triangle $xyz$, and
adding a new edge between $x$ and one of the two neighbors of $v_1$ in $B_1$
yields a graph $G'$ in ${\cal G}$ with $C(G')\geq C(G)+\frac{1}{n}$,
contradicting the choice of $G$.
Hence, since $G$ is not a tree, 
we may assume that $G$ has exactly one endblock $B$ that is a diamond, and all other blocks of $G$ are $K_2$s.
In this case $C(G)=\frac{8}{3n}$.
Since ${\cal G}$ contains a graph $G'$ with two endblocks that are triangles, we obtain $C(G')\geq \frac{14}{3n}>C(G)$,
contradicting the choice of $G$.
$\Box$

\medskip

\noindent Recall that the {\it paw} is the unique graph with degree sequence $1,2,2,3$.

\begin{claim}\label{claim7}
$|{\cal I}|\leq 1$.
\end{claim}
{\it Proof of Claim \ref{claim7}:}
Suppose, for a contradiction,  
that $G$ has two inner triangles $T_1$ and $T_2$.

If $T_1,T_2\in {\cal I}_3$, then 
contracting both triangles to single vertices $u_1$ and $u_2$,
adding a new paw $P$, and then
replacing one of the edges incident with $u_1$, say $u_1v$,
with the two new edges $xu_1$ and $xv$,
where $x$ is the vertex of degree $1$ in $P$,
yields a graph $G'$ in ${\cal G}$ with $C(G')=C(G)+\frac{1}{3n}$,
contradicting the choice of $G$.

If $T_1\in {\cal I}_2$ and $T_2\in {\cal I}_3$, then 
contracting both triangles to single vertices $u_1$ and $u_2$,
adding a new diamond $D$, and
adding the new edge $xu_1$,
where $x$ is a vertex of degree $2$ in $D$,
yields a graph $G'$ in ${\cal G}$ with $C(G')=C(G)$, $m(G')=m(G)$, 
and less inner triangles than $G$,
contradicting the choice of $G$.

If $T_1,T_2\in {\cal I}_2$, then
replacing $T_1$ with an edge between the two vertices, say $a$ and $b$, outside of $T_1$ that have neighbors in $T_1$,
adding a new triangles $xyz$, and 
adding the new edge $xu_2$,
where $u_2$ is the vertex in $T_2$ of degree $2$ in $G$, 
yields a graph $G'$ in ${\cal G}$ with $C(G')=C(G)$, $m(G')=m(G)$, 
and less inner triangles than $G$,
contradicting the choice of $G$.
Note that, 
in this last construction, 
the vertices $a$ and $b$ are non-adjacent in $G$ by Claim \ref{claim3}, and 
two triangles that contributed to ${\cal I}_2$ are replaced by one that contributes to ${\cal I}_3$.
$\Box$

\medskip

\noindent For $t(G)=(|{\cal D}|,|{\cal I}_2|,|{\cal I}_3|)$, the above claims imply 
$t(G)\in \{ 
(0,0,0),
(1,0,0),
(2,0,0),
(0,0,1),
(0,1,0)\}.$
For each of these cases, we can determine $C(G)$ exactly.
Let $k$ be the number of vertices of $G$ that do not lie in a triangle.

If $t(G)=(i,0,0)$ for some $i\in\{ 0,1,2\}$,
then $G$ arises from a tree of order $2k+2$ 
with $k$ vertices of degree $3$ and $k+2$ endvertices, 
by replacing $k+2-i$ endvertices with triangles and $i$ endvertices with diamonds,
and we obtain 
$n=4k+6+i$, and
$$C(G)
=\frac{1}{n}\left(\frac{7}{3}(k+2-i)+\frac{8}{3}i\right)
=\frac{7}{12}+\frac{14-3i}{12n}.$$
If $t(G)=(0,0,1)$,
then $G$ arises from a tree of order $2k+4$
with $k+1$ vertices of degree $3$ and $k+3$ endvertices,
by replacing all endvertices as well as one internal vertex with triangles,
and we obtain 
$n=4k+12$, and
$$C(G)
=\frac{1}{n}\left(\frac{7}{3}(k+3)+1\right)
=\frac{7}{12}+\frac{1}{n}.$$
If $t(G)=(0,1,0)$,
then $G$ arises from a tree of order $2k+3$
with $k$ vertices of degree $3$, one vertex of degree $2$, and $k+2$ endvertices,
by replacing all endvertices as well as the vertex of degree $2$ 
with triangles,
and we obtain 
$n=4k+9$, and
$$C(G)
=\frac{1}{n}\left(\frac{7}{3}(k+2)+\frac{5}{3}\right)
=\frac{7}{12}+\frac{13}{12n}.$$
Considering the different parities of $n$ modulo $4$,
the desired result follows.
$\Box$

\medskip

\noindent Our next goal is the characterization of all extremal graphs for (\ref{e1}).
Therefore, let ${\cal B}_0$ be the set of all connected subcubic graphs $G$ of order at least $6$ 
such that every block of $G$ is $K_2$, $K_3$, or a diamond, and every block of $G$ that is a diamond is an endblock of $G$. 
Let the {\it type} $t(G)$ of a graph $G$ in ${\cal B}_0$ be the $3$-tuple $(d,i_2,i_3)$, where 
$d$ is the number of blocks of $G$ that are diamonds, and, for $j$ in $\{ 2,3\}$,
$i_j$ is the number of blocks of $G$ that are triangles that contain exactly $j$ vertices of degree $3$ in $G$.

Let 
$${\cal B}
=
\Big\{ G\in {\cal B}_0:t(G)\in 
\{ 
(0,0,0),
(1,0,0),
(0,1,0),
(0,0,1),
(0,1,1),
(0,2,0),
(0,3,0)\}\Big\}.$$

\begin{theorem}\label{theorem3}
If $G$ is a connected subcubic graph of order at least $6$, 
then $G$ satisfies (\ref{e1}) with equality if and only if $G\in {\cal B}$.
\end{theorem}
{\it Proof:} First, let $G\in {\cal B}$.
For $t(G)\in \{ (0,0,0),(1,0,0),(0,1,0),(0,0,1)\}$, 
we verified at the end of the proof of Theorem \ref{theorem2} that (\ref{e1}) holds with equality. 
For $t(G)\in \{ (0,1,1),(0,2,0),(0,3,0)\}$, 
very similar simple calculations imply the same.

Now, let $G$ satisfy (\ref{e1}) with equality. 
Similarly as above, let ${\cal G}$ be the set of all connected subcubic graphs of order $n(G)$.
We consider the claims from the proof of Theorem \ref{theorem2}.
Clearly, Claim \ref{claim1} still holds.
Suppose, for a contradiction, that Claim \ref{claim2} fails, 
that is, $G$ has some edge $uv$ that lies in some cycle but in no triangle.
Iteratively removing from $G$ first the edge $uv$, and then further edges that lie in cycles but not in triangles as long as possible
yields a graph $G'$ in ${\cal G}$ with $C(G')\geq C(G)$ such that every edge of $G'$ that lies in some cycle also lies in some triangle.
Now, the argument from the proof of Claim \ref{claim3} applies, and, hence, $G'\in {\cal B}_0\subseteq {\cal G}$.
By (\ref{e1}) for $G'$, we obtain $C(G')=C(G)$, 
which, as observed in the proof of Claim \ref{claim2}, implies that $u$ and $v$ are vertices of $G'$ 
that do not lie in some triangle and have degree at most $2$ in $G'$.
Arguing as in the proof of Claim \ref{claim6}, 
we obtain the existence of a graph $G''$ in ${\cal G}$ with $C(G'')>C(G')=C(G)$,
which contradicts (\ref{e1}) for $G''$.
Hence, Claim \ref{claim2} holds.
Considering their respective proofs, 
it follows that also Claims \ref{claim3}, \ref{claim4}, \ref{claim5}, and \ref{claim6} hold.   
Let the type $t(G)$ of $G$ be $(d,i_2,i_3)$.
By Claim \ref{claim5}, $d\leq 2$.
If $d=2$, then, by Claim \ref{claim4}, $i_2=i_3=0$, 
and the calculation at the end of the proof of Theorem \ref{theorem2} implies the contradiction that (\ref{e1}) does not hold with equality.
Hence, $d\in \{ 0,1\}$, and, if $d=1$, then, by Claim \ref{claim4}, $i_2=i_3=0$.
If $i_3\geq 2$, then arguing as in the proof of Claim \ref{claim7} 
yields the existence of a graph $G'$ in ${\cal G}$ with $C(G')>C(G)$,
which is a contradiction.
Hence, $i_3\leq 1$.
If $i_2\geq 4$, then arguing as in the proof of Claim \ref{claim7} 
yields the existence of a graph $G'$ in ${\cal G}$ with $C(G')=C(G)$ that is of type $(d',i_2',i_3')=(d,i_2-4,i_3+2)$.
Now, as observed above, $i_3'\geq 2$ implies the existence of a graph $G''$ in ${\cal G}$ with $C(G'')>C(G')=C(G)$,
which is a contradiction.
Hence, $i_2\leq 3$.
Finally, 
if $i_2\geq 2$ and $i_3=1$, then arguing as in the proof of Claim \ref{claim7} 
yields the existence of a graph $G'$ in ${\cal G}$ with $C(G')=C(G)$ that is of type $(d',i_2',i_3')=(d,i_2-2,i_3+1)$.
Again, as observed above, $i_3'\geq 2$ implies a contradiction.
Hence, $i_3=1$ implies $i_2\leq 1$.
Altogether, it follows that $G\in {\cal B}$,
which completes the proof. $\Box$

\medskip

\noindent Our final result shows that adding a single edge can increase the clustering coefficient of a graph from $0$ to almost $1$,
which means that it is a rather sensitive parameter.

\begin{theorem}\label{theorem4}
If $G$ is a graph of order $n$ at least $3$, and $u$ and $v$ are non-adjacent vertices in $G$, then
$$C(G+uv)\leq C(G)+\left(1-\frac{2}{n}+\frac{4}{n(n-1)}\right)$$
with equality if and only of $G$ is $K_{2,n-2}$, and $u$ and $v$ are of degree $n-2$ in $G$.
\end{theorem}
{\it Proof:} Since the statement is trivial if $d_G(u)=0$ or $d_G(v)=0$, we may assume that neither $u$ nor $v$ are isolated in $G$.
Let $G'=G+uv$.

If $w$ is a common neighbor of $u$ and $v$ in $G$, then the degree of $w$ in $G'$ equals the degree of $w$ in $G$ 
but $G'[N_{G'}(w)]$ contains exactly one edge more than $G[N_G(w)]$, which implies
\begin{align*}
C_w(G') - C_{w}(G) &=\frac{1}{{d_G(w)\choose 2}}\leq 1
\end{align*}
with equality if and only if $d_G(w) = 2$.

If $d_G(u) \geq 2$, then
\begin{align*}
C_u(G') - C_u(G) &= \frac{m(G'[N_{G'}(u)])}{{d_{G'}(u)\choose 2}} -\frac{m(G[N_G(u)])}{{d_G(u)\choose 2}}\\
&= \frac{m(G[N_{G}(u)]) + |N_G(u) \cap N_G(v)|}{{d_G(u)+1\choose 2}} -\frac{m(G[N_G(u)])}{{d_G(u)\choose 2}}\\
&= \frac{2(d_G(u) -1)|N_G(u) \cap N_G(v)|- 4 m(G[N_G(u)])}{(d_{G}(u)+1)d_{G}(u)(d_{G}(u)-1)}\\
&\leq \frac{2(d_G(u) -1)d_G(u)- 4 m(G[N_G(u)])}{(d_{G}(u)+1)d_{G}(u)(d_{G}(u)-1)}\\
&\leq \frac{2}{d_{G}(u)+1}
\end{align*}
with equality if and only if $N_G(u)\subseteq N_G(v)$ and $N_G(u)$ is independent.

If $d_G(u)=1$, then 
$$C_u(G') - C_u(G)\leq 1=\frac{2}{d_G(u)+1}$$
with equality if and only if $N_G(u)\subseteq N_G(v)$ and $N_G(u)$ is independent.

We obtain that 	
\begin{align*}
n(C(G')-C(G)) &= (C_u(G')-C_u(G))+(C_v(G')-C_v(G))+\sum_{w\in N_G(u) \cap N_G(v)}(C_w(G')-C_w(G))\\
& \leq \frac{2}{d_G(u)+1}+\frac{2}{d_G(v)+1}+|N_G(u)\cap N_G(v)|
\end{align*}
with equality if and only if 
$N_G(u)=N_G(v)$ is independent and every vertex in $N_G(u)$ has degree $2$ in $G$,
that is, $G$ is $K_{2,n-2}$, and $u$ and $v$ are of degree $n-2$ in $G$.
Since $d:=\min\{ d_G(u),d_G(v)\}\leq n-2$ and the function $f:\mathbb{N}\to\mathbb{R}:d\mapsto \frac{4}{d+1}+d$ is strictly increasing, we conclude
\begin{align*}
n(C(G')-C(G)) & \leq \frac{2}{d_G(u)+1}+\frac{2}{d_G(v)+1}+|N_G(u)\cap N_G(v)|\\
& \leq \frac{4}{d+1}+d\\
& \leq \frac{4}{n-1}+n-2,
\end{align*}
which completes the proof. $\Box$

\end{document}